\documentclass[11pt]{article}
\usepackage{graphicx}
\usepackage{amsmath}
\usepackage{amssymb}
\usepackage{epstopdf}
\usepackage[small,nohug,heads=vee]{diagrams}
\diagramstyle[labelstyle=\scriptstyle]

\def\eproof{$\Box$ \medskip}

\newcommand{\B}{\mathcal B}
\newcommand{\G}{\mathfrak g}

\newcommand{\R}{\mathbb R}
\newcommand{\SL}{{\mathsf{SL}}(n,\R)}
\newcommand{\PSL}{{\mathsf{PSL}}(n,\R)}
\newcommand{\sln}{{\mathfrak sl}(n,\R)}
\newcommand{\gl}{{\mathfrak gl}(n,\R)}
\newcommand{\Rep}{\mathcal R}

\newcommand{\rpn}{{\mathbb{RP}}^{n-1}}
\newcommand{\Hp}{{\mathbb H}^2}

\newcommand{\oa}{\overline{\alpha}}
\newcommand{\ob}{\overline{\beta}}

\newtheorem{theorem}{Theorem}
\newtheorem{lemma}{Lemma}

\title{The Poisson Bracket of Length functions in the Hitchin Component}
\author{Martin Bridgeman}
\begin{document}
\maketitle

\begin{abstract}
Wolpert's cosine formula on Teichm\"uller space gives the Weil-Petersson Poisson bracket $\{\l_\alpha, l_\beta\}$ for geodesic length functions $l_\alpha,l_\beta$ of closed curves $\alpha,\beta$  as the sum of the cosines of the angle of intersection of the associated geodesics. This was recently generalized to Hitchin representations by Labourie. In this paper, we give a short proof of this generalization using Goldman's formula for the Poisson bracket on representation varieties of surface groups into reductive Lie groups.
\end{abstract}

\section{Introduction}
Let $S$ be a closed oriented surface of genus $g \geq 2$. In \cite{hitchin}, Hitchin considered the space  
$$\Rep_n(S) = Hom^{red}(\pi_1(S),\PSL)/\PSL$$ 
of conjugacy classes of reducible representations of $\pi_1(S)$ into $\PSL$.  The space $\Rep_n(S)$ has the structure of algebraic variety.

In \cite{goldman-components}, Goldman showed that for $n=2$, $\Rep_2(S)$ has $4g-3$ components, two of which are the Teichm\"uller components $T(S), T(\overline{S})$ corresponding to the conformal structures on $S$ and its complex conjugate $\overline{S}$ respectively. The space $\Rep_n(S)$ has a natural symplectic structure $\omega$, called the {\em Goldman symplectic form},  discovered  by Goldman (see \cite{goldman-symplectic}). This generalized the symplectic form discovered by Atiyah-Bott for the case of representations into the group $U(n)$ (see \cite{AB}). For $n=2$ the  form $\omega$ restricts on $\Rep_n(S)$ to  (an integer multiple of) the well-known Weil-Petersson symplectic form  $\omega_{wp}$ on $T(S)$. 

The symplectic form $\omega$ on $\Rep_n(S)$  defines a dual Poisson structure on $\Rep_n(S)$ given by $\{f,g\} = \omega(Hf,Hg)$ where $Hf, Hg$ are the Hamiltonian vector fields  with respect to $\omega$ of the smooth functions $f,g:\Rep_n(S) \rightarrow \R$.

Given $\alpha$ a  homotopy class of a non-trivial closed curve on $S$, we have the associated length function $l_\alpha:T(S)\rightarrow \R$ which assigns  the length of the geodesic representative of $\alpha$ in the associated hyperbolic structure. In \cite{wolpert-fn}, Wolpert  showed that  for the Weil-Petersson symplectic form, then $Hl_\alpha = -t_\alpha$ where $t_\alpha$ is the twist vector field obtained by dehn twist about $\alpha$ a simple non-trivial closed curve. Wolpert further proved the the following {\em cosine formula} for the Poisson bracket of length functions.

\begin{theorem}{(Wolpert, \cite{wolpert-cosine}})
Let  $\{.,.\}_{wp}$ be the Poisson bracket on Teichm\"uller space $T(S)$ given by the Weil-Petersson symplectic form. Let $\alpha, \beta$ be homotopy classes of closed oriented curves in $S$ with unique closed geodesic representatives $\oa,\ob$ in $X \in T(S)$. Then
$$\{l_{\alpha},l_{\beta}\}_{wp}(X) = \sum_{p \in \oa \cap \ob} \cos \theta_p$$
where $\theta_p$ is the angle of intersection of $\oa, \ob$ at $p$ measured  from $\oa$ to $\ob$ counterclockwise.
\end{theorem}

As part of his proof of the Nielsen realization conjecture (see \cite{kerckhoff-nielsen}), Kerkhoff also derived the above formula for the case when the curves are measured laminations. 

In the recent preprint, {\em Goldman algebra, opers and the swapping algebra}, Labourie generalizes the above formula for Hitchin representations (see \cite[Theorem 6.1.2]{labourie-anosov}). In this note, we give another proof of this generalization using Goldman's formula for the Poisson bracket of invariant functions (see \cite{goldman-bracket}).

A representation $\rho:\pi_1(S) \rightarrow \PSL$ is {\em Hitchin} if there exists a Teichm\"uller representation $\rho_0:\pi_1(S) \rightarrow PSL(2,\R)$ such that  $\rho = \tau_n \circ \rho_0$ where $\tau_n:\mathsf{PSL}(2,\R) \rightarrow \PSL$ is the irreducible representation. As Teichm\"uller space is connected, Hitchin representations correspond to  (at most two) connected components of $R_n(S)$ given by the images of $T(S), T(\overline{S})$ under $\tau_n$. Thus for $n = 2$ the Hitchin components are exactly the Teichm\"uller components $T(S), T(\overline{S})$.  Hitchin proved the following;

\begin{theorem}{(Hitchin, \cite{hitchin})}
Each Hitchin component  is homeomorphic to $\R^{|\chi(S)|(n^2-1)}$. If $n$ is even there are exactly two Hitchin components and if $n$ is odd, there is exactly one.
\end{theorem}

Using techniques from the dynamics of Anosov flows, Labourie showed the following;

\begin{theorem}{(Labourie, \cite{labourie-anosov})}
If $\rho$ is a Hitchin representation then $\rho$ is discrete faithful and for every $g \neq e$, $\rho(g)$  is diagonalizable over $\R$ with eigenvalues distinct $\lambda_1(g),\ldots,\lambda_n(g)$ satisfying
$$|\lambda_1(g)| > |\lambda_2(g)| >\ldots > |\lambda_n(g)|.$$
\end{theorem}

Thus given $\alpha$ a homotopy class or closed oriented curve in $S$, we therefore have functions $l^i_{\alpha}:H_n(S) \rightarrow \R$ given by 
$$l^i_{\alpha}([\rho]) = \log |\lambda_i(\rho(\alpha))|.$$

 In \cite{labourie-cross}, Labourie introduced the following {\em cross-ratio} on quadruples of lines and planes.  We let $\rpn$ be the space of lines in $\R^n$ (considered as  non-zero vectors in $\R^n$ up to multiplication by $\R^*$), and ${\rpn}^*$ the space of planes (considered as the space of  non-zero linear functionals on $\R^n$ up to multiplication by $\R^*$). 
 
 The cross-ratio is given by the map $b: {\rpn}\times{\rpn}^*\times\rpn\times {\rpn}^*$ 
$$b(x,y,z,w) = \frac{<y'|z'><w'|x'>}{<y'|x'><w'|z'>}$$
where $x' \in x, y' \in y, z' \in z, w' \in w$ are any choice of non-zero elements. By linearity  $b$ is well defined as the above formula is independent of the choices made. The cross-ratio $b$ is obviously only defined when the quadruple $(x,y,z,w)$ is in general position.

Given  $A$ a matrix with eigenvalues having distinct absolute values, we define
$\xi^i(A) \in \rpn$ to be the $i$-th eigenspace, and $\theta^i(A) \in {\rpn}^*$ to be the plane spanned by $\{\xi^j\}_{j\neq i}$. We let $\xi(A) = (\xi^1(A), \xi^2(A),\ldots,\xi^n(A))$ and $\theta(A) = (\theta^1(A),\theta^2(A),\ldots,\theta^n(A))$. We define
$$b^{ij}(A,B) =  b(\xi^i(A), \theta^i(A),\xi^j(B),\theta^j(B)).$$
If $\rho:\pi_1(S) \rightarrow \PSL$ is a Hitchin representation, and $\alpha, \beta \in \pi_1(S)$ then we define 
$$b_\rho^{ij}(\alpha,\beta) = b^{ij}(\rho(\alpha),\rho(\beta)).$$

In \cite{labourie-swap}, Labourie gives the following generalization of Wolpert's cosine formula.
\begin{theorem}{(Labourie, \cite{labourie-swap})}
Let $\alpha, \beta$ be homotopy classes of closed oriented curves in $S$ represented by  immersed curves $\oa,\ob$ in $S$ which are in general position, then
$$\{l^i_\alpha,l^j_\beta\}([\rho]) = \sum_{p \in \oa\cap\ob} \epsilon(p,\oa,\ob) \left(b^{ij}_{\rho_p}(\oa_p,\ob_p)- \frac{1}{n}\right).$$
\label{gwolpert}
\end{theorem}

We will give an elementary proof of this theorem.

We note that for $n = 2$ there is a single cross-ratio $b$ and for $A, B \in {\mathsf {PSL}}(2,\R)$, $b(A,B) = \cos^2(\phi_p/2)$ where $\phi_p$ is the angle of intersection between the positive rays of the associated geodesics in $\alpha,\beta$ in  $\Hp$  a the point of intersection $p = \alpha \cap \beta$.  
Thus 
$$b(A,B) - \frac{1}{2} = \frac{1}{2}(2\cos^2(\phi_p/2)-1) = \frac{1}{2}\cos(\phi_p).$$

The angle $\theta_p < \pi$  is  the counterclockwise angle between $\alpha,\beta$ at their intersection point.  Thus if $0 < \phi_p < \pi$,  $p$ is positively oriented then $ \phi_p = \theta_p$ and if $\pi <\phi_p < 2\pi$ then $p$ is negatively oriented and $\theta_p = \phi_p - \pi$. Thus the above formula for $n =2$ is
$$\{l^1_\alpha,l^1_\beta\}([\rho]) = \sum_{p \in \oa\cap\ob} \epsilon(p,\oa,\ob) \left(\frac{1}{2}\cos(\phi_p)\right) =  \frac{1}{2}\sum_{p \in \oa\cap\ob}  \cos(\theta_p).$$

For $g \in \SL(2,\R)$  we have $\lambda_1(g) = e^{l(g)/2}$ where $l(g)$ is the hyperbolic translation of $g$. Therefore it follows that if $l_\gamma$ is the length function for closed curve $\gamma$ then $ l_\gamma = 2l^1_\gamma $. Also the classical Weil-Petersson symplectic form $\omega_{wp}$ satisfies  $\omega =  2\omega_{wp}$ (see \cite{goldman-bracket}). Therefore we recover Wolpert's cosine formula for the Weil-Petersson Poisson structure
$$\{l_\alpha, l_\beta\}_{wp} = \sum_{p \in \alpha' \cap \beta'} \cos \theta_p.$$

\section{Background}
We now describe the background on Goldman's formula for the Poisson bracket of invariant functions.
Let $G$ be a reductive matrix group and consider the non-degenerate symmetric form $\B:\G \times \G \rightarrow \R$ given by $\B(X,Y) = Tr(XY)$. An invariant function for $G$ is a smooth function $f:G\rightarrow \R$ which is conjugacy invariant. In particular $f = Tr$ is an invariant function. Given $f$ there is a natural function $F:G \rightarrow \G$ given by
$$\B(F(A), X)  = \frac{d}{dt}f(exp(tX)A)\qquad \mbox{for all } X\in\G$$
Thus $F(A)$ is dual to $R^*_A(df(A)) \in \G^*$ under the isomorphism  $\hat{\B}: \G\rightarrow \G^*$ given by $\hat{\B}(X)(Y) = \B(X,Y)$.

Let $S$ be a closed oriented surface of genus $g \geq 2$ and $\pi = \pi_1(S, p)$ for some $p \in S$. We consider the space $Hom(\pi,G)/G$ of representations $\rho:\pi \rightarrow G$ up to conjugacy and let $\Rep(S,G)$ be the space of smooth points of $Hom(\pi,G)/G$.  If $\alpha$ is a non-trivial homotopy class of closed oriented  curve in $S$ then  $\alpha$ defines a conjugacy class  in $\pi$. If $f$ is an invariant function for $G$ then we can define $f_\alpha:\R(S,G) \rightarrow \R$ by
$$f_\alpha([\rho]) = f(\rho(\alpha'))$$
where $\alpha' \in \alpha$. 

The tangent space at $[\rho] \in \Rep(S,G)$ can be identified with the group cohomology $H^1(\pi, \G_{Ad\circ \rho})$. Using $\B$ to pair coefficients, we use the cup-product and cap-product  for group cohomology to define the map
$$H^1(\pi, \G_{Ad\circ \rho})\times H^1(\pi, \G_{Ad\circ \rho}) \xrightarrow{\B(.\cup.)} H^2(\pi,\R) \xrightarrow{\cap[\pi]} H_0(\pi,\R) = \R$$
 This map  defines  the  Goldman symplectic form $\omega$ on $\Rep(S,G)$ (see \cite{goldman-symplectic}). Specifically we have
$$\omega_{[\rho]}(X,Y) = \B(X\cup Y)\cap [\pi].$$
Given a smooth function $f:\Rep(S,G)\rightarrow \R$ the {\em Hamiltonian vector field} of $f$ is the vector field $Hf$ defined by
$\omega(Hf,Y) = df(Y)$. For $f,g$ two smooth functions the associated {\em Poisson bracket}  on smooth functions is the pairing $\{.,.\}:C^\infty(\Rep(S,G),\R)\times C^\infty(\Rep(S,G),\R)   \rightarrow C^\infty(\Rep(S,G)),\R)$ given by
$$\{f,g\} ([\rho])= \omega_{[\rho]}(Hf,Hg).$$

Given $\alpha$ an oriented curve in $S$, if $p \in \alpha$, we let $\alpha_p$ be the oriented curve given by traversing $\alpha$ starting at $p$.  If $\alpha,\beta$ are two  oriented closed curves, then $\alpha, \beta$ are in general position if their intersections are transverse. If $\alpha,\beta$ are in general position, then for $p \in \alpha\cap \beta$ we define $\epsilon(p,\alpha,\beta) = \pm 1$ given by if the orientation of the point of intersection agrees or not with the orientation of the surface.

Also  for $[\rho] \in \Rep(S,G)$ we  let $\rho_p:\pi_1(S,p) \rightarrow G$ be a representation defined by change of base point of $\rho$. This is well-defined up to conjugacy.

Goldman gave the following description of the Poisson bracket  for invariant functions.

\begin{theorem}{(Goldman, \cite{goldman-bracket})}
Let $f,f': G \rightarrow \R$ be invariant functions for $G$ with associated functions $F,F':G\rightarrow \G$. Let $\alpha, \beta$ be homotopy classes of closed oriented curves represented by  immersed curves $\oa,\ob$ in $S$ which are in general position. Then
$$\{f_{\alpha},f'_{\beta}\} [\rho] = \sum_{p\in \oa\cap \ob} \epsilon(p,\oa,\ob)\B(F(\rho_p(\oa_p)), F'(\rho_p(\ob_p))$$
\label{poisson}
\end{theorem}

\section{Length Functions}
As Hitchin representations can be lifted to representations into $\SL(n,\R)$ (see \cite{labourie-anosov}), we can restrict to representations into $\SL(n,\R)$.  We define the hyperbolic elements  $Hyp \subseteq \SL$ to be the open subset  of diagonalizable matrices over $\R$ with eigenvalues having distinct absolute values. For $A \in Hyp$, $A$ has eigenvalues $\lambda_1(A),\ldots, \lambda_n(A)$  with $|\lambda_1(A)| > |\lambda_2(A)| > \ldots>|\lambda_n(A)|$.
We define the functions $l^i:Hyp \rightarrow \R$ by letting $L^i(A) = \log|\lambda_i(A)|$. 
We define the function $L^i: Hyp \rightarrow \sln$ by
$$\B(L^i(A),X) = \frac{d}{dt} l^i(\exp(tX)A).$$

\subsection{Eigenvalue Perturbation}

We now consider perturbation of eigenvalues in the space of hyperbolic matrices. Given $A \in Hyp$ let $p_i(A):\R^n \rightarrow \R^n$ be projection onto the $i$-th eigenspace, parallel to the other eigenvectors.  

\begin{lemma}
The length function $l^i:Hyp \rightarrow \R$ satisfies
$$dl^i_A(X) = \frac{1}{\lambda^i(A)}Tr(p_i(A).X).$$
\end{lemma}

{\bf Proof:}
We let $A = A_0$ and denote the eigenvalues and eigenvectors   of $A$ by $\lambda^i, x^i$. We further let  $\dot{A} = \dot{A}_0$.
 We have  $A_t$ has eigenvalues $\lambda^i_t$ and  unit length eigenvector $x^i_t$. We have
$$A_t.x^i_t = \lambda^i_t.x^i_t.$$
Differentiating we get
$$ \dot{A}x^i + A\dot{x}^i = \lambda^i.\dot{x}^i+\dot{\lambda}^i.x^i $$
We let $p_i(A)$ be linear projection onto the $i$-th eigenspace of $A$ parallel to the other eigenspaces of $A$.
We apply to the above equation.
$$ p_i(A)\dot{A}x^i + p_i(A)A\dot{x}^i = \lambda^i.p_i(A)\dot{x}^i+\dot{\lambda}^i.x^i $$
As $p_i(A)A =\lambda^i.p_i(A)$
we have $p_i(A)A\dot{x}^i = \lambda^i.p_i(A)\dot{x}^i$ so after cancellation we get
$$\dot{\lambda}^i.x^i = p_i(A).\dot{A}.x^i.$$
Therefore we have
$$\dot{\lambda}^i = tr(p_i(A).\dot{A}).$$
As $l^i(X) = \log|\lambda^i(X)|$  on $Hyp$ we have
$$dl^i = \frac{d\lambda^i}{\lambda^i}.$$
Therefore
$$dl^i_A(X) = \frac{1}{\lambda^i(A)}Tr(p_i(A).X)$$
\eproof

We now use the above lemma to calculate $L^i$.

\begin{lemma}
$$L^i(A) = p_i(A) - \frac{1}{n}I.$$
\end{lemma}

{\bf Proof:}
By the above
$$\B(L^i(A),X) = \frac{d}{dt} l^i(\exp (tX)A) = dl^i_A(XA) = \frac{1}{\lambda^i(A)}Tr(p_i(A).XA).$$
By definition $A.p_i(A) = \lambda^ip_i(A)$. Therefore
$$\B(L^i(A),X) = \frac{1}{\lambda^i(A)}Tr(A.p_i(A).X) =  \frac{1}{\lambda^i(A)}Tr(\lambda^i.p_i(A).X) = Tr(p_i(A).X).$$
We let $\overline{\B}:\gl\times \gl \rightarrow \R$ given by $\B(X,Y) = Tr(XY)$. Then $\overline{\B}$ is non-degenerate and restricts to $\B$ on $\G$.
We let $P:\gl \rightarrow \sln$ be orthogonal projection with respect to  $\overline{\B}$. Then given $A \in \gl$, then for all $X \in \G$
$$\overline{\B}(A,X) = \overline{\B}(P(A),X) = \B(P(A),X).$$
Therefore we have
$$\B(L^i(A),X) = Tr(p_i(A).X) = \overline{\B}(p_i(A),X) =  \B(P(p_i(A)),X).$$
As $\B$ is non-degenerate on $\sln$, we have
$$L^i(A) = P(p_i(A))$$
The projection map $P:\gl \rightarrow \sln$ is given by
$$P(A) = A - \frac{1}{n}Tr(A).I$$
Therefore as $p_i(A)$ is projection onto a 1-dimensional eigenspace, $Tr(p_i(A)) = 1$ and we have
$$L^i(A) = p_i(A) - \frac{1}{n}Tr(p_i(A)).I = p_i(A) - \frac{1}{n}.I$$
\eproof

\subsection{Poission bracket}
We now use Goldman's formula to give an alternative proof of Labourie's generalization of the cosine formula.

{\bf Theorem \ref{gwolpert}} {\em (Labourie, \cite{labourie-swap})}
{\em
$$\{l^i_\alpha,l^j_\beta\}([\rho]) = \sum_{p \in \oa \cap \ob} \epsilon(p,\oa,\ob) \left(b^{ij}_{\rho_p}(\oa_p, \ob_p) - \frac{1}{n}\right).$$
}

{\bf Proof:}
From the above we have
$$\B(L^i(A), L^j(B)) = Tr\left(\left(p_i(A)-\frac{1}{n}.I\right).\left(p_j(B)-\frac{1}{n}.I\right)\right) $$
As $Tr(p_i(A)) = Tr(p_j(B)) = 1$
$$\B(L^i(A), L^j(B)) = Tr(p_i(A)p_j(B)) - \frac{1}{n}$$
Now applying Goldman's formula from Theorem \ref{poisson} we get
$$\{l^i_\alpha,l^j_\beta\}([\rho]) = \sum_{p \in \oa\cap \ob} \epsilon(p,\oa,\ob) \B(L^i(\rho(\oa_p)),L^j(\rho(\ob_p))) $$
$$= \sum_{p \in \oa\cap \ob} \epsilon(p,\oa,\ob) \left(Tr(p_i(\rho(\oa_p))p_j(\rho(\ob_p)) - \frac{1}{n}\right).$$

For any $X \in Hyp$ and let $\xi(X), \theta(X)$ be the $n$-tuples of eigenspaces and dual planes.
We let $A,B \in Hyp$ and we choose non-zero elements $a^i_+ \in \xi^i(A), a^i_- \in \theta^i(A), b^j_+ \in \xi^j(B), b^j_- \in \theta^j(B)$. Then
$$p_i(A)(v)  = \frac{<a^i_-| v>}{<a^i_-| a^i_+>}a^i_+ \qquad p_j(B)(v)  = \frac{<b^j_-| v>}{<b^j_-| b^j_+>}b^j_+.$$
Similarly for $B \in Hyp$ with $b^+_i, b^-_i$. Then if $A, B \in Hyp$ we have
$$p_i(A)p_j(B) v =  \frac{<a^i_-| b^j_+>}{<a^i_-| a^i_+>}\frac{<b^j_-| v>}{<b^j_-| b^j_+>}$$
Thus 
$$Tr(p_i(A)p_j(B)) =   \frac{<a^i_-| b^j_+><b^j_-| a^i_+>}{<a^i_-| a^i_+><b^j_-| b^j_+>} = b(\xi^i(A),\theta^i(A),\xi^j(B),\theta^j(B)) = b^{ij}(A,B).$$
Therefore the Poisson bracket is
$$\{l^i_\alpha,l^j_\beta\}([\rho]) = \sum_{p \in \oa \cap \ob} \epsilon(p,\oa,\ob) \left(b^{ij}_{\rho_p}(\oa_p, \ob_p) - \frac{1}{n}\right).$$
\eproof

\end{document}